\documentclass[leqno]{article}
\usepackage{amsmath,amsfonts,amsthm,amssymb,indentfirst}

\newtheorem{theorem}{Theorem}
\newtheorem{lemma}[theorem]{Lemma}
\newtheorem{definition}[theorem]{Definition}

\newcommand{\E}{\mathrm{End}}

\begin{document}

\title{Generating Endomorphism Rings of Infinite Direct Sums and Products of Modules}
\author{Zachary Mesyan}
\date{May 12, 2004}

\maketitle

\section{Introduction}
In~\cite{DM&PN} Macpherson and Neumann show that for any infinite
set $\Omega$ the full permutation group $\mathrm{Sym}(\Omega)$ is
not the union of a chain of $\leq |\Omega|$ proper subgroups,
while in~\cite{Bergman} Bergman shows that if $U$ is a generating
set for $\mathrm{Sym}(\Omega)$ as a group, then there exists a
positive integer $n$ such that every element of
$\mathrm{Sym}(\Omega)$ can be written as a group word of length
$\leq n$ in elements of $U$. We will modify the arguments used in
these two papers to prove analogous statements for endomorphism
rings of infinite direct sums and products of isomorphic modules.
Namely, we will show that if $R$ is a unital associative ring, $M$
is a left $R$-module, $\Omega$ is an infinite set, and $N$ is
either $\bigoplus_{i \in \Omega} M$ or $\prod_{i \in \Omega} M$,
then the ring $\E_R (N)$ is not the union of a chain of $\leq
|\Omega|$ proper subrings, and also that given a generating set
$U$ for $\E_R (N)$ as a ring, there exists a positive integer $n$
such that every element of $\E_R (N)$ is represented by a ring
word of length at most $n$ in elements of $U$. (The notion of ring
word will be made precise below.)

\section{Moieties}
First we will fix some notation and prove two somewhat technical
lemmas.  For the rest of this section $R$ will denote a ring, $M$
a nonzero left $R$-module, and $\Omega$ an infinite set. $N$ will
denote either $\bigoplus_{i \in \Omega} M$ or $\prod_{i \in
\Omega} M$ (all the arguments will work under either
interpretation), and $E$ will denote $\E_R (N)$. Endomorphisms
will be written on the right of their arguments. Also, given a
subset $\Sigma \subseteq \Omega$, we will write $M^{\Sigma}$ for
the $R$-submodule of $N$ corresponding to $\Sigma$, so in
particular, $N=M^{\Sigma} \oplus M^{\Omega\backslash\Sigma}$. If
$\Sigma \subseteq \Omega$ and $U \subseteq E$, then let
$U_{\{\Sigma\}} = \{f \in U : M^{\Sigma}f \subseteq M^{\Sigma},\
M^{\Omega \backslash \Sigma} f \subseteq M^{\Omega \backslash
\Sigma}\}$ and $U_{[\Sigma]} = \{f \in U : M^{\Sigma}f \subseteq
M^{\Sigma},\ M^{\Omega \backslash \Sigma}f = \{0\}\}$. We will
also say that a subset $\Sigma \subseteq \Omega$ is \emph{full}
with respect to $U \subseteq E$ if the set of endomorphisms of
$M^{\Sigma}$ induced by members of $U_{\{\Sigma\}}$ is all of
$\E_R (M^{\Sigma})$. Finally, $\Sigma \subseteq \Omega$ is called
a \emph{moiety} if $|\Sigma| = |\Omega| = |\Omega \backslash
\Sigma|$.

We will now prove analogs of Lemmas 3 and 4 of~\cite{Bergman}.

\begin{lemma}\label{x,y}
Suppose that a subset $U \subseteq E$ has a full moiety
 $\Sigma \subseteq \Omega$. Then there exist elements $x, y \in
 E$ such that $E = yUy + yUyx + xyUy + xyUyx$.
\end{lemma}

\begin{proof}
Let $\pi_{\Sigma}$ denote the projection from $N$ to $M^{\Sigma}$
along $M^{\Omega\backslash\Sigma}$ and $\pi_{\Omega \backslash
\Sigma}$ denote the projection from $N$ to $M^{\Omega \backslash
\Sigma}$ along $M^{\Sigma}$. Then for any $f \in E$ we can write
$f = \pi_{\Sigma}f\pi_{\Sigma} + \pi_{\Sigma}f\pi_{\Omega
\backslash \Sigma} + \pi_{\Omega \backslash \Sigma}f\pi_{\Sigma} +
\pi_{\Omega \backslash \Sigma}f\pi_{\Omega \backslash \Sigma}$. We
also note that $\pi_{\Sigma}U\pi_{\Sigma} = E_{[\Sigma]}$, since
$\Sigma$ is full with respect to $U$.

Now $|\Sigma| = |\Omega \backslash \Sigma|$, as $\Sigma$ is a
moiety, so there is an automorphism $x \in E$ of order 2 such that
$M^{\Sigma}x = M^{\Omega \backslash \Sigma}$ and $M^{\Omega
\backslash \Sigma}x = M^{\Sigma}$. Then
$\pi_{\Sigma}f\pi_{\Sigma}$, $\pi_{\Sigma}f\pi_{\Omega \backslash
\Sigma}x$, $x\pi_{\Omega \backslash \Sigma}f\pi_{\Sigma}$,
$x\pi_{\Omega \backslash \Sigma}f\pi_{\Omega \backslash \Sigma}x
\in E_{[\Sigma]}$.  Hence, $\pi_{\Sigma}f\pi_{\Omega \backslash
\Sigma} = \pi_{\Sigma}f\pi_{\Omega \backslash \Sigma}x^2 \in
E_{[\Sigma]}x$, $\pi_{\Omega \backslash \Sigma}f\pi_{\Sigma} =
x^2\pi_{\Omega \backslash \Sigma}f\pi_{\Sigma} \in xE_{[\Sigma]}$,
and $\pi_{\Omega \backslash \Sigma}f\pi_{\Omega \backslash \Sigma}
= x^2\pi_{\Omega \backslash \Sigma}f\pi_{\Omega \backslash
\Sigma}x^2 \in xE_{[\Sigma]}x$.  So $f \in E_{[\Sigma]} +
E_{[\Sigma]}x + xE_{[\Sigma]} + xE_{[\Sigma]}x =
\pi_{\Sigma}U\pi_{\Sigma} + \pi_{\Sigma}U\pi_{\Sigma}x +
x\pi_{\Sigma}U\pi_{\Sigma} + x\pi_{\Sigma}U\pi_{\Sigma}x$.
\end{proof}

The proof of the following lemma is set-theoretic in nature, so
aside from a few minor adjustments, we present it here the way it
appears in~\cite{Bergman}.

\begin{lemma}\label{diagonal}
Let $(U_i)_{i \in I}$ be any family of subsets of $E$ such that
$\bigcup_{i \in I} U_i = E$ and $|I| \leq |\Omega|$. Then $\Omega$
contains a full moiety with respect to some $U_i$.
\end{lemma}

\begin{proof}
Since $|\Omega|$ is infinite and $|I| \leq |\Omega|$, we can write
$\Omega$ as a union of disjoint moieties $\Sigma_i$, $i \in I$.
Suppose that there are no full moieties with respect to $U_i$ for
any $i \in I$.  Then in particular, $\Sigma_i$ is not full with
respect to $U_i$ for any $i \in I$.  So for every $i \in I$ there
exists an endomorphism $f_i \in \E_R (M^{\Sigma_i})$ which is not
the restriction to $M^{\Sigma_i}$ of any member of
$(U_i)_{\{\Sigma_i\}}$. Now if we take $f \in E$ to be the
endomorphism whose restriction to each $M^{\Sigma_i}$ is $f_i$,
then $f$ is not in $U_i$ for any $i \in I$, contradicting
$\bigcup_{i \in I} U_i = E$.
\end{proof}

\section{Generating Sets}
We are now ready to prove our main results.

\begin{theorem}\label{chain}
Let $R$ be a ring, $M$ a nonzero left $R$-module, $\Omega$ an
infinite set, and $E = \E_R (\bigoplus_{i \in \Omega} M)$ $($or
$\E_R (\prod_{i \in \Omega} M))$. Suppose that $(R_i)_{i \in I}$
is a chain of subrings of $E$ such that $\bigcup_{i \in I} R_i =
E$ and $|I| \leq |\Omega|$. Then $E = R_i$ for some $i \in I$.
\end{theorem}

\begin{proof}
By the preceding lemma, $\Omega$ contains a full moiety with
respect to some $R_i$. Thus, Lemma~\ref{x,y} implies that $E =
\langle R_i \cup \{x, y \} \rangle$ for some $x, y \in E$. But, by
the hypotheses on $(R_i)_{i \in I}$, $R_i \cup \{x, y \} \subseteq
R_j$ for some $j \in I$, and hence $E = \langle R_i \cup \{x, y \}
\rangle \subseteq R_j$, since $R_j$ is a subring.
\end{proof}

\begin{definition}\label{word}
Let $R$ be a ring and $U$ a subset of $R$.  We will say that $r
\in R$ is represented by a \emph{ring word of length 1} in
elements of $U$ if $r \in U \cup \{0, 1, -1\}$, and, recursively,
that $r \in R$ is represented by a \emph{ring word of length} $n$
in elements of $U$ if $r = p + q$ or $r = pq$ for some elements
$p, q \in R$ which can be represented by ring words of lengths
$m_1$ and $m_2$ respectively, with  $n = m_1 + m_2$.
\end{definition}

\begin{theorem}\label{gen}
Let $R$ be a ring, $M$ a nonzero left $R$-module, $\Omega$ an
infinite set, and $E = \E_R (\bigoplus_{i \in \Omega} M)$ $($or
$\E_R (\prod_{i \in \Omega} M))$. If $U$ is a generating set for
$E$ as a ring, then there exists a positive integer $n$ such that
every element of $E$ is represented by a ring word of length at
most $n$ in elements of $U$.
\end{theorem}

\begin{proof}
For $i = 1, 2, 3, \dots,$ let $U_i$ be the set of elements of $R$
that can be expressed as ring words in elements of $U$ of length
$\leq i$. Then $\bigcup_{i = 1}^{\infty} U_i = E$, since $U$ is a
generating set. Since $\aleph_0 \leq |\Omega|$,
Lemma~\ref{diagonal} implies that $\Omega$ contains a full moiety
with respect to some $U_i$. By Lemma~\ref{x,y}, there exist $x, y
\in E$ such that $E = yU_iy + yU_iyx + xyU_iy + xyU_iyx$. Let $k
\geq i$ be such that $U_k$ contains $U_i \cup \{x, y \}$.  We then
have $E = U_k^3 + U_k^4 + U_k^4 + U_k^5$.

Now, for any positive integer $m$, $U_k^m$ consists of ring words
of length $\leq mk$ in elements of $U$. Thus $E = U_k^3 + U_k^4 +
U_k^4 + U_k^5$ consists of ring words of length $\leq 16k$ in
elements of $U$.
\end{proof}

\section{Acknowledgment}
I would like to thank Professor Bergman for bringing this topic to
my attention and for his suggestions and comments throughout the
preparation of this note.

Department of Mathematics, University of California, Berkeley, CA,
94720

Email: zak@math.berkeley.edu


\begin{thebibliography}{00}
\bibitem{Bergman} George M. Bergman,
{\it Generating infinite symmetric groups,} preprint, 8pp.,
http://math.berkeley.edu/{$\!\sim$}gbergman/papers/Sym\_Omega:1,
arXiv:\linebreak[0]math.GR/0401304.

\bibitem{DM&PN} H. D. Macpherson and Peter M. Neumann,
{\it Subgroups of infinite symmetric groups,} J. London Math. Soc.
(2) {\bf 42} (1990) 64--84. ~MR~{\bf 92d}:20006.~
\end{thebibliography}
\end{document}